\documentclass[12pt]{article}
\usepackage{amsmath,amsfonts,amssymb,amsthm}
\usepackage{bbm,mathrsfs,bm,mathtools}
\title{Measure transport via pseudo-cones}
\author{Rolf Schneider}
\date{}
\newcommand{\Sn}{{\mathbb S}^{n-1}}

\newcommand{\R}{{\mathbb R}}

\newcommand{\N}{{\mathbb N}}

\newcommand{\D}{{\rm d}}

\newtheorem{theorem}{Theorem}
\newtheorem{lemma}{Lemma}

\begin{document}
\maketitle

\begin{abstract}
{For the solution of the Gauss image problem for pseudo-cones, which can be considered as a measure transport problem for certain measures on the sphere, we give a new proof, using a special case of Kantorovich duality.}\\[2mm]
\end{abstract}

\section{Introduction}\label{sec1}

Let $C\subset\R^n$ be a closed convex cone, pointed (that is, not containing a line) and with interior points. A $C$-pseudo-cone is a nonempty closed convex set $K\subset C$ not containing the origin and satisfying $K+C =K$. This implies that $C$ is the recession cone of $K$ and that $\lambda K\subseteq K$ for $\lambda\ge 1$. The set of all $C$-pseudo-cones in $\R^n$ is denoted by $ps(C)$. 

Denoting by $\Sn$ the unit sphere of $\R^n$, we define the sets $\Omega_C:= \Sn\cap{\rm int}\,C$ and $\Omega_{C^\circ}:= \Sn\cap{\rm int}\,C^\circ$, where $C^\circ:=\{x\in\R^n: \langle x,y\rangle\le 0\,\forall y\in C\}$ is the dual cone of $C$. (We denote by ${\rm int}$ the interior, by ${\rm cl}$ the closure, by ${\rm bd}$ the boundary, by $\langle\cdot\,,\cdot\rangle$ the scalar product, and by $o$ the origin of $\R^n$.) 

Let $K\in ps(C)$. Then $\rho_K(v):= \min\{r>0: rv\in K\}$, $v\in\Omega_C$, defines the {\em radial function} $\rho_K:\Omega_C\to(0,\infty)$ of $K$ (the set $\{r>0: rv\in K\}$ is not empty for $v\in\Omega_C$, as follows from $K=K+C$). Obviously, $\rho_K(v)v\in {\rm bd}\,K$.

Let $\omega_K$ be the set of all $u\in{\rm cl}\,\Omega_{C^\circ}$ with the property that $u$ is a normal vector (by which we always mean an outer unit normal vector) at more than one point of $K$. Then $\omega_K$ has Hausdorff dimension at most $n-2$ (as follows from \cite[Thm. 2.2.5]{Sch14} and polarity). The {\em reverse radial Gauss map} $\alpha_K^*:\Omega_{C^\circ}\setminus \omega_K\to \Omega_C$ is defined by letting $\alpha_K^*(u)$ be the unique vector $v\in\Omega_C$ for which $u$ is attained as a normal vector of $K$ at $\rho_K(v)v$. The mapping $\alpha_K^*$ is continuous and hence measurable.

For a topological space $X$, we denote by $P(X)$ the set of Borel probability measures on $X$. The main subject of this note is the following result.

\begin{theorem}\label{T1.1}
Let $\mu\in P(\Omega_{C^\circ})$ and $\nu\in P(\Omega_C)$ and suppose that $\mu$ is zero on sets of Hausdorff dimension at most $n-2$. Then there exists a $C$-pseudo-cone $K\in ps(C)$ such that $(\alpha_K^*)_\sharp\mu=\nu$.
\end{theorem}

As usual, $(\alpha_K^*)_\sharp\mu=\nu$ means that  $\mu((\alpha_K^*)^{-1}(\eta)) = \nu(\eta)$ for each Borel set $\eta\subseteq \Omega_C$, that is, $\nu$ is the image measure of $\mu$ under the reverse radial Gauss map of $K$.

Theorem \ref{T1.1} was proved in \cite{Sch25a}. It was motivated by the introduction and investigation of the Gauss image problem for convex bodies in \cite{BLYZZ20}. The proof there used a variational argument, as was the case in \cite{Sch25a}. A treatment of the Gauss image problem for convex bodies by mass transportation methods was presented by Bertrand \cite{Ber23}. The purpose of the following is to prove also Theorem \ref{T1.1} by such methods. The proof is much shorter than that in \cite{Sch25a}, though less elementary. We do not try to carry over the argumentation of \cite{Ber23} to pseudo-cones, but apply the Kantorovich duality in a similar way as done by R\"uschendorf and Rachev \cite{RR90}. 

As in \cite{Sch25a}, $\nu$ is a measure on $\Omega_C$ and not on its closure, although the latter case would also be of interest. The reason for this restriction in  the present proof is that for securing the assumption in Theorem \ref{T3.1} we have first to consider measures $\nu$ on a compact subset of $\Omega_C$.

Theorem \ref{T1.1} is proved in Section \ref{sec3}, while Section \ref{sec2} contains some preparations.

\section{Preparations involving pseudo-cones}\label{sec2}

In this section, we collect and define some notions which will be essential for the proof of Theorem \ref{T1.1}. Let $K\in ps(C)$.

The support function of $K$ is defined by
$$ h_K(u)= \sup_{y\in K} \langle u,y\rangle\quad \mbox{for }u\in\Omega_{C^\circ},$$
from which it follows that
$$ h_K(u)= \sup_{v\in\Omega_C}\langle u,v\rangle\rho_K(v).$$
Writing $\overline h_K=-h_K$, this is equivalent to
\begin{equation}\label{M1}
\overline h_K(u)=\inf_{v\in\Omega_C}|\langle u,v\rangle|\rho_K(v).
\end{equation}

The copolar set of $K$ is defined by
$$ K^*= \{x\in\R^n: \langle x,y\rangle\le -1\mbox{ for all }y\in K\}.$$
It satisfies $K^*\in ps(C^\circ)$, and we have $K^{**}=K$ (proved in \cite[Thm. 3.9]{XLL23}). Moreover,
\begin{equation}\label{M2}
\rho_K(v)=\frac{1}{\overline h_{K^*}(v)}= \frac{1}{\inf_{u\in\Omega_{C^\circ}} |\langle u,v\rangle|\rho_{K^*}(u)},
\end{equation}
where the first equality is shown, for example, in \cite{XLL23} and the second equality follows from (\ref{M1}).

Let $f: \Omega_C\to (0,\infty)$ be a function that is bounded away from $0$ (that is, satisfies $f\ge a$ for some $a>0$). We define the {\em convexification} of $f$ by
$$ 
\langle f\rangle:= \bigcap\{K\in ps(C): f(v)v\in K\,\forall v\in\Omega_C\}.
$$
Since $f$ is bounded away from $0$, there exists $K\in ps(C)$ with $f(v)v\in K$ for all $v\in\Omega_C$ (for example, the intersection of $C$ with a suitable closed halfspace). Clearly, $\langle f\rangle\in ps(C)$,
\begin{equation}\label{M3} 
\rho_{\langle f\rangle} \le f
\end{equation}
and
\begin{equation}\label{M4} 
\langle \rho_K\rangle =K \quad\mbox{for } K\in ps(C).
\end{equation}

For the support function of $\langle f\rangle$ we get, by (\ref{M1}) and (\ref{M3}),
$$ \overline h_{\langle f\rangle}(u)=\inf_{v\in\Omega_C}|\langle u,v\rangle|\rho_{\langle f\rangle}(v) \le \inf_{v\in\Omega_C}|\langle u,v\rangle|f(v).$$
Here strict inequality would contradict the definition of the convexification, hence we have
\begin{equation}\label{M5}
\overline h_{\langle f\rangle}(u)=\inf_{v\in\Omega_C}|\langle u,v\rangle|f(v)\quad\mbox{for }u\in\Omega_{C^\circ}.
\end{equation}

Similarly, of course, for $g:\Omega_{C^\circ}\to(0,\infty)$, bounded away from $0$, we define
$$ \langle g\rangle:= \bigcap \{K\in ps(C^\circ): g(u)u\in K\,\forall u\in\Omega_{C^\circ}\},$$
which is a $C^\circ$-pseudo-cone.

For a function $f:\Omega_C\to (0,\infty)$, bounded away from zero, we define the {\em pseudo-conjugate} function $f^*:\Omega_{C^\circ}\to (0,\infty)$ by
$$ f^*(u):= \sup_{v\in\Omega_C} \frac{1}{|\langle u,v\rangle|f(v)},\quad u\in\Omega_{C^\circ}.$$
For fixed $u\in\Omega_{C^\circ}$, the product $|\langle u,v\rangle| f(v)$ is bounded away from $0$, hence $f^*(u)<\infty$. Trivially, $f^*$ is also bounded away from zero, so that $f^{**}$ can be defined. Of course, for functions $g:\Omega_{C^\circ}\to (0,\infty)$, we define $g^*$ correspondingly.

It follows from (\ref{M2}) and (\ref{M5}) that
$$ \rho_{\langle f \rangle^*}(u) = \frac{1}{\overline h_{\langle f\rangle}(u)} = \frac{1}{\inf_{v\in\Omega_C} |\langle u,v\rangle| f(v)}= f^*(u),$$
thus
\begin{equation}\label{M6}
f^*=\rho_{\langle f\rangle^*}.
\end{equation}
Since the pseudo-conjugate function of $f$ is the radial function of the copolar set of the convexification of $f$, the function $f^*$ is always the radial function of a $C^\circ$-pseudo-cone. It follows from (\ref{M6}) and (\ref{M4}) that 
\begin{equation}\label{M6a}
(\rho_K)^*=\rho_{\langle\rho_K\rangle^*} =\rho_{K^*} 
\end{equation}
for $K\in ps(C)$.

Again, let $f:\Omega_C\to(0,\infty)$ be a function that is bounded away from $0$. Using (\ref{M6}), (\ref{M6a}) and (\ref{M3}), we get
$$ f^{**}=(\rho_{\langle f\rangle^*})^*= \rho_{\langle f\rangle^{**}}= \rho_{\langle f\rangle} \le f,$$ thus
\begin{equation}\label{M7}
f^{**} \le f.
\end{equation}
Here equality holds if and only if $f$ is the radial function of a $C$-pseudo-cone.

From the definition of $f^*$ it follows that
$$ f^*(u) \ge \frac{1}{|\langle u,v\rangle|f(v)},$$
hence
\begin{equation}\label{M8}
f(v)f^*(u) \ge \frac{1}{|\langle u,v\rangle|} \quad\mbox{for }(v,u)\in\Omega_C\times\Omega_{C^\circ}.
\end{equation}
Suppose that equality holds in (\ref{M8}). Then, by (\ref{M6}), (\ref{M2}) and (\ref{M3}),
$$
f(v) = \frac{1}{|\langle u,v\rangle|f^*(u)} = \frac{1}{|\langle u,v\rangle|\rho_{\langle f\rangle^*}(u)} \le \frac{1}{\overline h_{\langle f\rangle^*}(v)} = \rho_{\langle f\rangle}(v)\le f(v).
$$
Since equality holds here, we have $ \rho_{\langle f\rangle}(v)=f(v)$ and $|\langle u,v\rangle|\rho_{\langle f \rangle^*}(u) =\overline h_{\langle f\rangle^*}(v)$. The latter means that $v$ is a normal vector of $\langle f\rangle^*$ at $\rho_{\langle f\rangle^*}(u)u$. Equivalently (as follows, e.g., from \cite[Lem. 6]{Sch24a}),
$u$ is a normal vector of $\langle f\rangle$ at $\rho_{\langle f\rangle}(v)v$. The converse also holds.

We recall that in \cite{Sch25b} we have defined the pseudo-subdifferential of $K$ by
$$ \partial^\bullet K:=\{(v,u)\in\Omega_C\times\Omega_{C^\circ}: \mbox{$u$ is a normal vector of $K$ at $\rho_K(v)v$}\}.$$
We can thus state the following lemma.

\begin{lemma}\label{L2.1}
Let $f:\Omega_C\to(0,\infty)$ be a function bounded away from zero. The inequality
\begin{equation}\label{M9} 
f(v)f^*(u)\ge\frac{1}{|\langle u,v\rangle|}\end{equation}
holds for all $(v,u)\in\Omega_C\times\Omega_{C^\circ}$. Equality in $(\ref{M9})$ for some $(v,u)\in\Omega_C\times\Omega_{C^\circ}$ holds if and only if $\rho_{\langle f\rangle}(v) =f(v)$ and $(v,u)\in \partial^\bullet \langle f\rangle$.
\end{lemma}

\section{Proof of Theorem \ref{T1.1}}\label{sec3}

We define the function $c:\Omega_{C^\circ}\times\Omega_C\to (0,\infty)$ by
\begin{equation}\label{3.0} 
c(u,v):= -\log|\langle u,v\rangle|,\quad (u,v)\in \Omega_{C^\circ}\times\Omega_C.
\end{equation}
This is the cost function found in \cite{Sch25b}, but with a negative sign, to make it nonnegative. 

Given $\mu\in P(\Omega_{C^\circ})$ and $\nu\in P(\Omega_C)$, we denote by $M(\mu,\nu)$ the set of all probability measures $\pi$ on $\Omega_{C^\circ}\times\Omega_C$ with marginals $\mu$ and $\nu$, that is, satisfying $\pi(\omega\times\Omega_C)=\mu(\omega)$ and $\pi(\Omega_{C^\circ}\times\eta)=\nu(\eta)$ for all Borel sets $\omega\subseteq\Omega_{C^\circ}$ and $\eta\subseteq\Omega_C$. Then we define, with $c$ given by (\ref{3.0}),
$$ S(\mu,\nu):= \sup\left\{\int_{\Omega_{C^\circ}\times\Omega_C} c(u,v)\,\D\pi(u,v): \pi\in M(\mu,\nu)\right\},$$
\begin{eqnarray*} 
&& I(\mu,\nu):= \\
&& \inf\left\{\int_{\Omega_{C^\circ}} h\,\D\mu+\int_{\Omega_C} g\,\D\nu: h\in L^1(\mu),\, g\in L^1(\nu),\, h(u)+g(v)\ge c(u,v)\,\forall u,v\right\}
\end{eqnarray*}
and make use of the following theorem.

\begin{theorem}\label{T3.1}
We have $S(\mu,\nu)=I(\mu,\nu)$. The $\sup$ defining $S(\mu,\nu)$ is an attained maximum. Suppose there exist functions $h\in L^1(\mu)$ and $g\in L^1(\nu)$ satisfying $h(u)+g(v) \ge c(u,v)$ for all $(u,v)\in \Omega_{C^\circ}\times\Omega_C$. Then the $\inf$ defining $I(\mu,\nu)$ is an attained minimum.
\end{theorem}

This theorem follows from Sections 2.2 and 2.3 of the book by Rachev and R\"uschendorf \cite{RR98}. More precisely, the definitions of $S(\mu,\nu)$ and $I(\mu,\nu)$ are special cases of Definition 2.2.2 in \cite{RR98}, with $n=2$ and $h=c$, where $c(u,v):= -\log|\langle u,v\rangle|$. In our case, the latter function is real-valued, positive, and continuous. It is not required in \cite{RR98} that it be bounded. Theorem \ref{T3.1} is then a combination of Theorems 2.3.1(c), 2.3.10 and 2.3.12 in \cite{RR98} (note that our assumption in Theorem \ref{T3.1} implies that $I(\mu,\nu)<\infty$). The sections 2.2 and 2.3 of \cite{RR98} are based on a general version of the Kantorovich duality proved by Kellerer \cite{Kel84}. This generalization was used before in a paper by R\"uschendorf and Rachev \cite{RR90}, after which our argument is modeled.

We assume first that $\nu$ is concentrated on a closed subset of $\Omega_C$. Then there is a constant $a$ with $0< a\le 1$ and $|\langle u,v\rangle|\ge a$ and hence $-\log|\langle u,v\rangle|\le|\log a|$ for all $(u,v)\in \Omega_{C^\circ}\times {\rm supp}\,\nu$ (where ${\rm supp}\,\nu$ denotes the support of the measure $\nu$). It follows that there exist $h\in L^1(\mu)$ and $g\in L^1(\nu)$ satisfying $h(u)+g(v) \ge c(u,v)$ for all $(u,v)\in \Omega_{C^\circ}\times\Omega_C$. In particular, $I(\mu,\nu)<\infty$.

By Theorem \ref{T3.1} there exist functions $h$ on $\Omega_{C^\circ}$ and $g$ on $\Omega_C$ such that $h\in L^1(\mu)$, $g\in L^1(\nu)$, 
$$ h(u)+g(v)\ge -\log|\langle u,v\rangle|\quad\mbox{and}\quad \int_{\Omega_{C^\circ}} h\,\D\mu+ \int_{\Omega_C} g\,\D \nu= I(\mu,\nu).$$
Let
$$ \varphi(u):=e^{h(u)},\qquad \psi(v):= e^{g(v)}.$$
Since $h$ and $g$ are bounded from below, the positive functions $\varphi$ and $\psi$ are bounded away from zero. In particular, we can define
$$ f:= \varphi^{**}.$$
We have
$$ \varphi(u)\psi(v)\ge e^{-\log|\langle u,v\rangle|}$$
and hence
$$\psi(v) \ge \frac{1}{|\langle u,v\rangle|\varphi(u)}$$
for all $u\in \Omega_{C^\circ}$ and therefore $\psi(v)\ge \varphi^*(v)$ for $v\in\Omega_C$. This gives
\begin{eqnarray*}
\log \varphi(u) +\log\psi(v)&\ge& \log\varphi(u) +\log \varphi^*(v)\\
&\ge& \log\varphi^{**}(u)+\log \varphi^*(v)\\
&=& \log f(u) +\log f^*(v),
\end{eqnarray*}
where we have used (\ref{M7}) and $\varphi^*=\varphi^{***}$. By (\ref{M9}) we have
$$ \log f(u) +\log f^*(u)\ge -\log|\langle u,v\rangle|.$$
Thus, the pair $(\log f,\log f^*)$ is also a solution of the minimum problem, hence we have
$$ \int_{\Omega_{C^\circ}} \log f\,\D \mu+\int_{\Omega_C} \log f^*\,\D\nu = I(\mu,\nu).$$
By Theorem \ref{T3.1}, there exists $\pi^*\in M(\mu,\nu)$ for which
$$ \int_{\Omega_{C^\circ}\times\Omega_C} -\log|\langle u,v\rangle|\,\D\pi^*(u,v) = I(\mu,\nu),$$
hence
$$ \int_{\Omega_{C^\circ}\times\Omega_C} \left[-\log|\langle u,v\rangle|-\log f(u)-\log f^*(v)\right]\,\D\pi^*(u,v) =0.$$
Since the integrand is non-positive, we deduce that it is zero almost everywhere, that is,
$$ f(u)f^*(v) =\frac{1}{|\langle u,v\rangle|}$$
for $\pi^*$-almost all $(u,v)$.

By its definition, $f^*$ is the radial function of some $C$-pseudo-cone $K$. By Lemma \ref{L2.1}, $(v,u)\in\partial^\bullet K$ for $\pi^*$-almost all $(u,v)\in\Omega_{C^\circ}\times\Omega_C$. In the following, we write $\partial^\circ K := \{(u,v): (v,u)\in\partial^\bullet K\}$.

Now suppose that $\mu$ is zero on sets of Hausdorff dimension at most $n-2$. Then $\omega_K$ has $\mu$-measure zero. For each $u\in \Omega_{C^\circ}\setminus \omega_K$ there is a unique pair $(u,v)\in\partial^\circ K$, and we have $\alpha_K^*(u)=v$. Let $\eta\subseteq\Omega_C$ be a Borel set. Since $\pi^*$ has marginal measures $\mu$ and $\nu$, we get
\begin{eqnarray*}
\nu(\eta) &=& \pi^*(\{(u,v)\in\Omega_{C^\circ}\times\Omega_C: u\in\Omega_{C^\circ},\,v\in\eta\})\\
&=& \pi^*(\{(u,v)\in\partial^\circ K:u\in\Omega_{C^\circ},\, v\in\eta\})\\
&=& \pi^*(\{(u,v)\in\partial^\circ K: \alpha_K^*(u)=v,\,v\in\eta\}),
\end{eqnarray*}
where we have used that $\pi^*((\Omega_{C^\circ}\times\Omega_C)\setminus \partial^\circ K)=0$ and 
$$ \pi^*(\{(u,v)\in\partial^\circ K: \alpha_K^*(u)\mbox{ is not defined at } u\})=0.$$
Moreover,
\begin{eqnarray*}
\mu((\alpha_K^*)^{-1}(\eta)) &=& \pi^*(\{(u,v)\in\Omega_{C^\circ}\times\Omega_C: u\notin\omega_K,\,\alpha_K^*(u)\in\eta,\,v\in\Omega_C\})\\
&=& \pi^*(\{(u,v)\in\partial^\circ K:\alpha_K^*(u)=v,\, v\in\eta\}).
\end{eqnarray*}
Hence, $\mu((\alpha_K^*)^{-1}(\eta))=\nu(\eta)$.

So far we have assumed that ${\rm supp}\,\nu$ is a compact subset of $\Omega_C$. Now let $\nu\in P(\Omega_C)$ be as in Theorem \ref{T1.1}. We choose a sequence $(\eta_j)_{j\in \N}$ of compact sets satisfying $\nu(\eta_1)>0$, $\eta_j\subset \eta_{j+1}$ for all $j\in\N$ and $\bigcup_{j\in\N} \eta_j =\Omega_C$. Then we set
$$ \nu_j(\eta) :=\frac{\nu(\eta\cap\eta_j)}{\nu(\eta_j)} \quad\mbox{for Borel sets } \eta\subseteq \Omega_C.$$
As shown above, to each $j\in\N$ there exists a $C$-pseudo-cone $K_j$ such that $(\alpha_{K_j}^*)_\sharp\mu=\nu_j$.

We recall the convergence of pseudo-cones. Let $(K_j)_{j\in\N}$ be a sequence in $ps(C)$. Denoting by $B^n$ the unit ball of $\R^n$ with center $o$, we say that $K_j\to K$ (a closed convex set) as $j\to\infty$ if there is some $t_0>0$ with $K_j\cap t_0B^n\not=\emptyset$ for $j\in\N$ and, for each $t\ge t_0$, the sequence $(K_j\cap tB^n)_{j\in\N}$ of convex bodies converges, in the usual Hausdorff sense, to $K\cap tB^n$. By Lemma 1 of \cite{Sch24a}, every sequence of $C$-pseudo-cones whose distances from $o$ are bounded and bounded away from zero has a subsequence that converges to a closed convex set, and it is easy to see that this set is a $C$-pseudo-cone.

The pseudo-cones $K_j$ with $(\alpha_{K_j}^*)_\sharp\mu=\nu_j$ found above can be assumed to have distance $1$ from the origin, since multiplication of a pseudo-cone $K$ by a positive factor does not change the function $\alpha_K^*$. Then the sequence $(K_j)_{j\in\N}$ has a subsequence converging to a $C$-pseudo-cone $K_0$, and after changing the notation we can assume that $K_j\to K_0$. 

If we define
$$ \omega_\Box:= \bigcup_{j\in\N_0} \omega_{K_j},$$
then all functions $\alpha_{K_j}^*$, $j\in\N_0$, are defined on $\Omega_{C^\circ}\setminus\omega_\Box$ and hence $\mu$-almost everywhere. 

We state that
\begin{equation}\label{3.1}
\alpha_{K_j}^*(u)\to \alpha_{K_0}^*(u)\quad\mbox{as }j\to\infty, \mbox{ for $\mu$-almost all $u$}.
\end{equation}
For the proof, let $u\in\Omega_{C^\circ}\setminus\omega_\Box$. For each $j\in\N_0$, there is a unique $x_j=\rho_{K_j}(v_j)v_j\in{\rm bd}\,K_j$ at which $u$ is attained as a normal vector of $K_j$. We choose a number $t_0$ with $x_0\in{\rm int}(t_0B^n)$ ($x_0\notin\partial C$, since $u\in\Omega_{C^\circ}\setminus \omega_{K_0}$). Then $K_j\cap t_0 B^n\not=\emptyset$ for all $j\in\N_0$, since $K_j$ has distance $1$ from the origin. Let $y_j=x_j$ if $x_j\in t_0B^n$; otherwise let $y_j\in {\rm bd}(K_j\cap t_0 B^n)$ be any point at which $u$ is attained as an outer normal vector of $K_j\cap t_0 B^n$; in the latter case, necessarily $y_j\in{\rm bd}(t_0 B^n)$. The sequence $(y_j)_{j\in\N}$ is bounded and hence has a subsequence converging to some $x$. The point $x$ is necessarily a boundary point of $K_0\cap t_0B^n$, and $u$ is a normal vector of $K_0\cap t_0B^n$ at this point. If $x\not= x_0$, then the segment with endpoints $x$ and $x_0$ is contained in ${\rm bd}\, K_0$, so that $u\in \omega_{K_0}$, a contradiction. Thus, $x=x_0$. This shows also that every convergent subsequence of $(y_j)_{j\in\N}$ converges to $x_0$, hence this sequence is convergent to $x_0$. Since $x_0\in{\rm int}(t_0B^n)$, the sequence $(x_j)_{j\in\N}$ converges to $x_0$. This implies that $(v_j)_{j\in\N}$ converges to $v_0$, which proves (\ref{3.1}).

Now choose $k\in\N$, and let $g:\Omega_C\to\R$ be a bounded continuous function. Since $(\alpha_{K_k}^*)_\sharp \mu=\nu_k$, the change of variables formula gives
$$ \int_{\Omega_{C^\circ}} g\circ\alpha_{K_k}^*\,\D\mu = \int_{\Omega_C} g\,\D\nu_k= \frac{1}{\nu(\eta_k)} \int_{\eta_k} g\,\D\nu.$$
Since (\ref{3.1}) holds, the dominated convergence theorem shows that, as $k\to\infty$, the left side converges to  $\int_{\Omega_{C^\circ}} g\circ\alpha_{K_0}^*\,\D\mu$, while the right side converges to $\int_{\Omega_C} g\,\D\nu$. We deduce that
$$ \int_{\Omega_C} g\,\D(\alpha_{K_0}^*)_\sharp\mu = \int_{\Omega_{C^\circ}} g\circ\alpha_{K_0}^*\,\D\mu  = \int_{\Omega_C} g\,\D\nu.$$
Since this holds for all bounded continuous functions $g$ on $\Omega_C$, we can conclude that $(\alpha_{K_0}^*)_\sharp \mu=\nu$.

\noindent Author's address:\\[2mm]
Rolf Schneider\\Mathematisches Institut, Albert--Ludwigs-Universit{\"a}t\\D-79104 Freiburg i.~Br., Germany\\E-mail: rolf.schneider@math.uni-freiburg.de

\end{document}